\documentclass[11pt]{article}
\usepackage{amsmath}
\usepackage{amssymb}
\usepackage{cite}
\title{On the Darboux transformations and sequences of $p$-orthogonal polynomials}
\begin{document}
\author{D. Barrios Rolan\'{\i}a, J.C. Garc\'{\i}a-Ardila, D. Manrique\\\\
Dpto. de Matem\'atica Aplicada a la Ingenier\'{\i}a Industrial\\
Dpto. de Inteligencia Artificial\\
Universidad Polit\'ecnica de Madrid }
\date{}
\maketitle
\newtheorem{lem}{Lemma}
\newtheorem{prop}{Proposition}
\newtheorem{defi}{Definition}
\newtheorem{coro}{Corollary}
\newtheorem{rem}{Remark}
\newtheorem{nota}{Note}
\newtheorem{teo}{Theorem}
\def\adots{\mathinner{\mkern2mu\raise1pt\hbox{.}\mkern2mu\raise4pt\hbox{.}\mkern2mu\raise7pt\hbox{.}\mkern1mu}}
\newcommand{\se}[1]{\medskip \begin{center}{\section{#1}} \end{center}
\medskip }
\begin{abstract}
For a fixed $p \in \mathbb{N}$, sequences of polynomials $\{P_n\}$, $n \in \mathbb{N}$, defined by a $(p+2)$-term recurrence relation are related to several topics in Approximation Theory. A $(p+2)$-banded matrix $J$  determines the coefficients of the recurrence relation of any of such sequences of polynomials. The connection between these polynomials and the concept of orthogonality has been already established through a $p$-dimension vector of functionals. This work goes further in this topic by analyzing the relation between such vectors for the set of sequences $\{P_n^{(j)}\}$, $n \in N$, associated with the Darboux transformations $J^{(j)}$, $j=1, ..., p,$ of a given $(p+2)$-banded matrix $J$. 
\end{abstract}
\section{Introduction}
For a fixed $p\in\mathbb{N}$ we consider a sequence of polynomials $\{P_n\},\,n\in \mathbb{N},$ defined by a ($p+2$)-term recurrence relation
\begin{equation}
\left.
\begin{array}{r}
P_{n+1}(z)+(a_{n,n}-z)P_{n}(z)+\displaystyle\sum_{j=1}^p a_{n,n-j}P_{n-j}(z)=0\,,\quad n\in \mathbb{N}\,,\\
\\
P_{-p}=\cdots =P_{-1}=0\,,\quad P_0\equiv 1\,.
\end{array}
\right\}
\label{1}
\end{equation}
These polynomials are related with several topics such as Hermite-Pad\'e approximants and vector continued fractions (\cite{2018,Kaliaguine}). In particular, this kind of polynomials plays an essential role in the study of some integrable systems (see for instance  \cite{enviado,primero,fullKostant}).

%If $v^{(0)}(z):=(P_0(z),P_1(z),\ldots )^T$, then (\ref{1}) can be rewritten as the following formal matricial product,
%\begin{equation*}
%\left(J-zI\right)v^{(0)}(z)=0\,,
%\end{equation*}
%where $I$ is the semi-infinite identity matrix. 

The relation between the concept of orthogonality and the sequences of polynomials verifying a ($p+2$)-term recurrence relation was established in \cite{iseghem-1987} in the following well-known result.

\begin{lem}\label{Maroni} 
With the above notation, the following statements are equivalent.
\begin{itemize}
\item [(i)] $\{P_{n}\},\,n\in \mathbb{N} , $ verify (\ref{1}) with $a_{n,n-p} \neq 0$  for all $n\in \mathbb{N}$.
\item [(ii)] There exists a vector of functionals 
$
\nu=\left(\nu_1,\,\ldots , \nu_p\right),
$
where each $\nu_r \in {\cal P}^{\prime} \,, \, r=1,\ldots, p\,,$ is defined on the space of polynomials ${\cal P}$ verifying
\begin{equation}
\left\{
\begin{array}{lccc}
\nu_r\left[z^kP_{n}(z) \right]=0\,,& k=0,1,\ldots \,,& kp+r \le n\,,\quad n\in \mathbb{N}\,, \\\\
\nu_r\left[z^kP_{kp+r-1}(z) \right]\neq 0\,,& k=0,1,\ldots 
\end{array}
\right.
\label{ortogonalidad}
\end{equation}
\end{itemize}
\end{lem}

In the sequel, we call {\em vector of $p$-orthogonality} to any vector of functionals $\nu=\left(\nu_1,\,\ldots, \nu_p\right)$  verifying (\ref{ortogonalidad}). In this case, we say that $\{P_n\}$ is a {\em sequence of $p$-orthogonal polynomials} with respect to $\nu$.

For a sequence of polynomials $\{P_{n}\}$ defined in (\ref{1}), the sequence of linear functionals $\{{\cal L}_n\}\,,\,n=0,1,\ldots, $ given by
\begin{equation}
{\cal L}_j\left[P_i\right]=\delta_{i,j}\,,\quad i,j=0,1,\ldots,
\label{dual}
\end{equation}
plays an relevant role in the study of the orthogonality. $\{{\cal L}_n\}$ is called {\em dual sequence }, and it is the unique sequence of functionals verifying (\ref{dual}).  It is easy to check that the $p$ first terms ${\cal L}_r\,,\,r=0,1\ldots , p-1\,,$ of the dual sequence verify the orthogonality conditions (\ref{ortogonalidad}). This fact proves the existence of some vector of $p$-orthogonality associated with each arbitrary sequence $\{P_n\}$ of polynomials. However, the uniqueness of a vector functional as in (\ref{ortogonalidad})  is not guaranteed. In fact, P. Maroni characterized in \cite{Maroni89} these vectors of $p$-orthogonality  as $\left(\nu_1,\,\ldots , \nu_p\right)$ such that 
\begin{equation}
\begin{array}{lll}
\nu_1& = & \lambda_{1,0}{\cal L}_0\\
\nu_2& = & \lambda_{2,0}{\cal L}_0+ \lambda_{2,1}{\cal L}_1\\
& \vdots &\\
\nu_p& = & \lambda_{p,0}{\cal L}_0+ \lambda_{p,1}{\cal L}_1+\cdots +\lambda_{p,p-1}{\cal L}_{p-1}
\end{array}
\label{escalera}
\end{equation}
being $\lambda_{i,j}\in \mathbb{C} $ and $\lambda_{i,i-1}\neq 0$ for $i, j+1\in\{1,\ldots , p\}$.

Associated with (\ref{1}), it is possible to define the $(p+2)$-banded matrix $J$ whose entries are the coefficients of the recurrence relation,
\begin{equation}\label{2}
J = \left( \begin{array}{cccccc}
%\phantom{0}
 a_{0,0} & 1 & &    \\
 a_{1,0} & a_{1,1} & 1 & \\
 \vdots      & \vdots    &  \ddots   & \ddots      \\
 a_{p,0} &a_{p,1} & \cdots & a_{p,p} &1 & \\
0& a_{p+1,1} & & \ddots& \ddots & \ddots\\
& 0 &\ddots  \\
&  &\ddots  \\
\end{array} \right)\,,
\end{equation}
where we assume $a_{j+p,j}\neq 0\,,\,j=0,1,\ldots\,$

The use of discrete Darboux transformations was proposed in \cite{Matveev1,Matveev2} with the focus in the application to the Toda lattices. In \cite{enviado,primero,fullKostant,nuevo,2018,Dani} this study was extended to ($p+2$)-banded matrices as (\ref{2}). In the present work we are concerned about finding relations between the vectors of $p$-orthogonality associated with the Darboux transformations of such matrices \eqref{2}. Thereby this paper complements the analysis that has been done in \cite{nuevo} for the Geronimus transformations. We include here the following summary, with the more relevant concepts, for an independent reading.

Let $C\in \mathbb{C}$ be such that the main determinants of the infinite matrix $J-CI$ verify 
\begin{equation}
\det \left(CI_n-J_n \right)\neq 0 \text{ for each } n \in \mathbb{N}
\label{C}
\end{equation}
Due to the well-known fact that 
$P_n(z)=\det \left(zI_n-J_n \right),
$
\eqref{C} is equivalent to $C\in \mathbb{C}$ is not a zero of the sequence $\{P_n\}$.  
(Here and in the sequel, given a semi-infinite matrix $A$,  we denote by $A_n,$   $n\in \mathbb{N}$, the finite matrix of order $n$ formed with the first $n$ rows and columns of $A$.)  In these conditions, there exist two lower and upper triangular matrices, $L$ and $U$  respectively, 
\begin{equation*}%\label{1-}
L= \left( \begin{array}{ccccccc}
1\\
 l_{1,1} & 1 & &    \\
 \vdots &\ddots &\ddots\\
 l_{p,1} & l_{p,2} & \dots & 1 && \\
0 & l_{p+1,2}  & \ddots & & \ddots \\
& 0 &\ddots &\ddots  &\ddots \\
\end{array} \right)\,,
\end{equation*}
and 
\begin{equation*}
U = \left( \begin{array}{cccccc}
\gamma_1 & 1  &    \\
&  \gamma_{p+2} & 1  \\
&  &\gamma_{2p+3} & \ddots  \\
   &&& \ddots\\
\end{array} \right),
\label{U}
\end{equation*}
such that 
\begin{equation}
J-CI=LU
\label{3}
\end{equation}
is the unique factorization of $J-CI$ in these conditions (see \cite{Gantmacher} for details). In the sequel we assume $C\in\mathbb{C}$ fixed.

In \cite{Dani} the following factorization of $L$ was given.

\begin{lem}{\cite[Theorem 1, pp. 118]{Dani}}
In the above conditions, if $l_{p+i,i+1}\ne 0$ for each $i=0,1,\ldots $, then there exist $p$ bi-diagonal matrices 
\begin{equation}
L^{(j)} = \left( \begin{array}{cccccc}
  1 & &    \\
 \gamma_{j+1} & 1 & \\
   & \gamma_{p+j+2}  & 1 \\
 &&\gamma_{2p+j+3} &\ddots\\
 &&&\ddots
\end{array} \right),\,j=1,2,\ldots,p\,,
\label{Darboux2}
\end{equation}
 with $\gamma_{(k-1)p+j+k}\neq 0$ for all $k\in \mathbb{N}$, verifying
\begin{equation}
L=L^{(1)}L^{(2)}\cdots L^{(p)}\,.
\label{(b)}
\end{equation}
Moreover, for each $j\in \{1,2,\ldots,p-1\}$ it is possible to choose certain set of $p-j$ elements $\gamma_{j+1}, \gamma_{p+j+2}, \ldots , \gamma_{(p-j)p}$ of $L^{(j)}$ such that the factorization (\ref{(b)}) is unique for the fixed set of $p(p-1)/2$ points
\begin{equation}
\left.
\begin{array}{ccccc}
\gamma_{2},&\gamma_{p+3} &\cdots & \cdots &\gamma_{(p-1)p},\\
\gamma_{3},&\gamma_{p+4} &\cdots & \gamma_{(p-2)p},\\
\vdots & \vdots &\adots \\
\gamma_{p-1},&\gamma_{2p},\\
\gamma_p.
\end{array}
\right\}
\label{puntos}
\end{equation}
\label{lema1}
\end{lem}

(\ref{3}) and (\ref{(b)}) provide the so called {\em Darboux factorization} of $J-CI$, defined as
\begin{equation}
J-CI=L^{(1)}L^{(2)}\cdots L^{(p)}U,
\label{5}
\end{equation}
where $L^{(1)},\,L^{(2)},\,\ldots L^{(p)}$ are bidiagonal matrices as in (\ref{Darboux2}) given in Lemma \ref{lema1} for certain set (\ref{puntos}) of $p(p-1)/2$ fixed entries. Each circular permutation of a Darboux factorization gives a new ($p+2$)-banded matrix,

\begin{equation}
J^{(j)}=CI+L^{(j+1)}L^{(j+2)}\cdots L^{(p)}U L^{(1)}L^{(2)}\cdots L^{(j)},\,j=1,2,\ldots, p\,.
\label{22}
\end{equation}

\begin{defi}
The permutations $J^{(j)}$ given in \eqref{22} are called Darboux transformation of $J-CI$. 
\end{defi}

As a consequence of the Darboux transformations, for each $j\in \{1,\ldots, p\}$  it is possible to define a new sequence of polynomials $\{P_n^{(j)}\}$ verifying a ($p+2$) recurrence relation,
\begin{equation}
\left.
\begin{array}{r}
P^{(j)}_{n+1}(z)+(a^{(j)}_{n,n}-z)P^{(j)}_{n}(z)+\displaystyle\sum_{s=1}^p a^{(j)}_{n,n-s}P^{(j)}_{n-s}(z)=0\,,\quad n\in \mathbb{N\,,}\\
\\
P^{(j)}_{-p}=\cdots =P^{(j)}_{-1}=0\,,\quad P^{(j)}_0\equiv 1\,.
\end{array}
\right\}
\label{11}
\end{equation}
 Another way to write (\ref{11}) is
 \begin{equation}
 \left(J^{(j)}-zI\right)v^{(j)}=0\,,
 \label{recuj}
 \end{equation}   
 where 
 $$
v^{(j)}(z)=\left( P^{(j)}_{0}(z), P^{(j)}_{1}(z),\ldots    \right)^T\,,\qquad j=0,1,\ldots, p\,.
$$
Here and in the following, we extend the notation taking $P^{(0)}_n=P_n\,,\,n\in \mathbb{N}$ and $J^{(0)}=J$.
   
Since Lemma \ref{Maroni}, in the above conditions there exists some vector of $p$-orthogonality
\begin{equation}
\nu^{(j)}=\left(\nu^{(j)}_1,\,\ldots , \nu^{(j)}_p\right)\,,\quad j=1,\ldots, p\,,
\label{funcional}
\end{equation}
such that the corresponding conditions of orthogonality like (\ref{ortogonalidad}) are verified.

In this work we analyze some relations between the vectors of $p$-orthogonality  
$
\nu=\left(\nu_1,\,\ldots , \nu_p\right)
$
associated to the polynomials $\{P_n\}$ and the vectors of $p$-orthogonality  (\ref{funcional}), $j=1,\ldots,p$, for a Darboux factorization of $J-CI$ under some conditions.  

\begin{defi}
	Here and in what follows,
	$(z-C)\nu_i ,\,i=1,\ldots, p,$ is a linear functional defined on the space $\mathcal{P}[z]$ of polynomials as
	$$
	(z-C)\nu_i[q]=\nu_i[(z-C)q]
	$$
	for each $q\in \mathcal{P}[z]$.
\end{defi}

We recall that it is possible to have several vectors of $p$-orthogonality associated with each sequence of polynomials verifying a $(p+2)$-term recurrence relation. If $\left(\nu_1,\,\ldots , \nu_p\right)$ verifies  (\ref{ortogonalidad}) then for $\left(c_1,\ldots, c_p\right)\in \mathbb{C}$ also $\left(c_1\nu_1,\,\ldots , c_p\nu_p\right)$ verifies (\ref{ortogonalidad}). In a more general procedure, the coefficients $\lambda_{i,j}$ of \eqref{escalera} and the determinants
\begin{equation*}
\Delta_m^{(j)}=\left|
\begin{array}{cccccccc}
\lambda_{j+2,0} & \lambda_{j+3,0}  &\cdots  &\cdots  & \lambda_{m,0}  &\cdots & \lambda_{j+m+1,0}\\
\vdots & \vdots & & & \vdots & \\
\lambda_{j+2,j+1} &\lambda_{j+3,j+1}  & \cdots  & \cdots  &\lambda_{m,1} &\cdots &\lambda_{j+m+1,j+1}\\
0 & \lambda_{j+3,j+2} &&& \vdots&& \vdots\\
\vdots & 0 &\ddots &&&  \\
\vdots & \vdots & \ddots &\ddots & \vdots&& \vdots  \\
0 & 0 & \cdots &0 & \lambda_{m,m-1} & \cdots & \lambda_{j+m+1,m-1}\\
\end{array}
\right|\,,
\label{determinantes}
\end{equation*}
$j=0,1,\ldots, p-1,\, m=1,2,\ldots,p-j-1,$
play an important role in this paper. Our main contribution is the following.

\begin{teo}
\label{principal}
Let 
$
\nu=\left(\nu_1,\,\ldots , \nu_p\right)
$
be a vector of $p$-orthogonality for $\{P_n\},\,n\in\mathbb{N},$ as in \eqref{escalera} such that 
\begin{equation}
\Delta_m^{(j)}\neq 0\,,\quad j=0,\ldots, p-1,\quad m=1,\ldots,p-j-1. 
\label{condition}
\end{equation}
Then there exists a Darboux factorization \eqref{5} of $J-CI$ such that, for each $j=1,2,\ldots, p$,
\begin{equation}
\nu^{(j)}=\left(\nu_{j+1},\ldots, \nu_p, (z-C)\nu_1,\ldots , (z-C)\nu_j\right)
\label{vectorj}
\end{equation}
is a vector of $p$-orthogonality for the sequence of polynomials $\{P_n^{(j)}\},\,n\in \mathbb{N},$ associated with the Darboux transformation $J^{(j)}$ of $J-CI$ given in \eqref{22} (where  we understand 
$\nu^{(p)}=\left((z-C)\nu_1,\ldots , (z-C)\nu_p\right)$ in \eqref{vectorj}).
\end{teo}

\begin{defi}
In the conditions of Theorem \ref{principal} we say that the vectors of $p$-orthogonality \eqref{vectorj} are the Darboux transformations of $\nu$.
\end{defi}

In Section~\ref{section2} some auxiliary results are established. In particular, some connections between the sequences $\{P_n^{(j)}\},\,j=0,\ldots, p,$ are given for each fixed Darboux factorization. Finally, Theorem \ref{principal} is proved in Section \ref{section3}.

\section{Orthogonality and  Darboux transformations}\label{section2}

In this section we assume that \eqref{5} is a fixed Darboux factorization of $J-CI$ corresponding with a given set of entries \eqref{puntos} of the matrices $L^{(1)}, \ldots, L^{(p)} $.

The following result establishes some relationships between the various sequences of polynomials associated with the Darboux transformations of $J-CI$. 
\begin{teo}
\label{teorema1}
We have the following relations between the sequences of polynomials $\{P^{(j)}_{n}\},\,j=0,1,\ldots ,p$.
\begin{eqnarray}
P^{(j+1)}_{n}(z) & = &\frac{P^{(j)}_{n+1}(z)+\displaystyle\sum_{s=0}^{p-1}g^{(j)}_{n+1,n-s+1}P^{(j)}_{n-s}(z)}{z-C}\,, \quad j=0,1,\ldots , p-1\,,\quad n\geq 0\,,\label{kernel1}\\\nonumber\\
P^{(p)}_{n}(z)&= & \frac{P_{n+1}(z)-\frac{P_{n+1}(C)}{P_{n}(C)}P_{n}(z)}{z-C}\,,  \quad  n\geq 0\,,
\label{kernel2}
\end{eqnarray}
where $g^{(j)}_{n+1,n-s+1}\in \mathbb{C}$ for $j=0,1,\ldots , p-1\,,\, n\geq 0$ and $s=0,\ldots, p-1.$
\end{teo}
\noindent
\underline{Proof}.- 
In  \cite{Dani} was proved that 
\begin{equation}
L^{(j+1)}L^{(j+2)}\cdots L^{(i)}v^{(i)}(z)=v^{(j)}(z),\quad 0\le j<i\le p\,,
\label{6}
\end{equation}
where the product of the triangular matrix $L^{(j+1)}L^{(j+2)}\cdots L^{(i)}$ times the vector $v^{(i)}(z)$ is  understanding in a formal sense.

Moreover, from this, (\ref{22}) and (\ref{recuj}),
%$$
%0= \left(J^{(j)}-zI\right)v^{(j)}(z)=(C-z)v^{(j)}(z)+L^{(j+1)}L^{(j+2)}\cdots L^{(p)}UL^{(1)}\cdots L^{(j)}v^{(j)}(z)\,.
%$$
%Then
\begin{equation}
L^{(j+1)}L^{(j+2)}\cdots L^{(p)}UL^{(1)}\cdots L^{(j)}v^{(j)}(z)=(z-C)v^{(j)}(z)\,,\quad j=0,1,\ldots, p\,,
\label{nueva}
\end{equation}
where we understand 
\begin{equation}
ULv^{(p)}(z)=(z-C)v^{(p)}(z)
\label{nuevap}
\end{equation}
 when $j=p$.
Replacing $j$ by $j+1$ in (\ref{nueva}),
$$
L^{(j+2)}L^{(j+3)}\cdots L^{(p)}UL^{(1)}\cdots L^{(j)}L^{(j+1)}v^{(j+1)}(z)=(z-C)v^{(j+1)}(z)\,,\quad j=0,1,\ldots, p-1\,,
$$
understanding $L^{(j+2)}L^{(j+3)}\cdots L^{(p)}=I$ when $j=p-1$.
From this and (\ref{6}) (for $i=j+1$),
\begin{equation}
L^{(j+2)}\cdots L^{(p)}UL^{(1)}\cdots L^{(j)}v^{(j)}(z)=(z-C)v^{(j+1)}(z)\,,\quad j=0,1,\ldots, p-1,
\label{nueva2}
\end{equation}
where 
$
L^{(2)}\cdots L^{(p)}Uv^{(0)}(z)=(z-C)v^{(1)}(z),
$ this is, $
L^{(1)}\cdots L^{(j)}=I
$ when $j=0$.

On the other hand, it is easy to check that, for each $k\in \mathbb{N}$,  the row $k$ of the infinite matrix $L^{(j+2)}\cdots L^{(p)}UL^{(1)}\cdots L^{(j)}$ is
$$
\left(g^{(j)}_{k,1},\,g^{(j)}_{k,2},\ldots , \,g^{(j)}_{k,k+1},\,0,\,\ldots \right),
$$
being $g^{(j)}_{k,k+1}=1$ and the entries $g^{(j)}_{k,s},\,s=1,\ldots,k,$ independent on $z$. Moreover, $g^{(j)}_{k,s}=0$ for $s\le k-p$  when $k>p$. In other words, $L^{(j+2)}\cdots L^{(p)}UL^{(1)}\cdots L^{(j)}$ is a $(p+1)$-banded Hessenberg matrix, where it is easy to see that $g^{(j)}_{k,k-p+1}\neq 0$ since $\gamma_r\neq 0\,\, (r\in \mathbb{N})$ .  Hence, taking into account (\ref{nueva2}),
$$
(z-C)P^{(j+1)}_{k-1}(z)=g^{(j)}_{k,k-p+1}P^{(j)}_{k-p}(z)+\cdots+g^{(j)}_{k,k+1}P^{(j)}_{k}(z)
$$
which drives to (\ref{kernel1}) when $k=n+1$. Note that, in this case, 
\begin{equation}
g^{(j)}_{n+1,n-p+2}\neq 0\,.
\label{asterisco}
\end{equation}

For $j=0$ and $i=p$, (\ref{6}) becomes 
\begin{equation*}
Lv^{(p)}(z)=L^{(1)}\cdots L^{(p)}v^{(p)}(z)=v^{(0)}(z)\,.
\label{0}
\end{equation*}

Therefore 
\begin{equation}
Uv^{(0)}(z)=(z-C)v^{(p)}(z)
\label{888}
\end{equation}
(see  (\ref{nuevap})) and, comparing the $(n+1)$-row in both sides of (\ref{888}),
\begin{equation}
(z-C)P^{(p)}_{n}(z)=\gamma_{n(p+1)+1}P_n(z)+P_{n+1}(z)
\label{88}
\end{equation}
The right hand side of (\ref{88}) is a polynomial with a root in $z=C$. Thus
\begin{equation}
\gamma_{n(p+1)+1}=-\frac{P_{n+1}(C) }{P_n(C)}\,.
\label{8888}
\end{equation}
From (\ref{88}) and (\ref{8888}) we arrive to (\ref{kernel2}). 
$\hfill \square$ 

\begin{rem}
(\ref{kernel1}) and (\ref{kernel2}) coincide in the classic case $p=1$. Both relations extend \cite[(7.3), pp. 35]{Chihara}, this is,
$$
P^{(1)}_{n}(z)= \frac{P_{n+1}(z)-\frac{P_{n+1}(C)}{P_{n}(C)}P_{n}(z)}{z-C},
$$
where the sequence of Kernel polynomials $\{P^{(1)}_{n}\}$ are defined in terms of  $\{P_{n}\}$. In this sense $\{P^{(j)}_{n}\}\,,\,j=1,\ldots, p,$ are extensions of this classical sequence of Kernel polynomials. 
\end{rem}

The above remark justifies the following definition. 

\begin{defi}
For each $j=1,\ldots, p$ the polynomials $\{P^{(j)}_{n}\}\,,n\in \mathbb{N},$ are called $j$-Kernel polynomials.
\end{defi}
For each sequence $\{P_n^{(j)}\}\,,j=0,1,\ldots ,p\,, $ we denote by $\{{\cal L}_n^{(j)}\}$ the corresponding dual sequence (taking ${\cal L}_n^{(0)}={\cal L}_n$). Equivalently to the behavior of the sequences of polynomials, the terms of the dual sequences are related to each other.  

\begin{lem}
With the above notation, for each $j=0,1\ldots, p-1$ and $n=0,1,\ldots $ we have
\begin{eqnarray}
{\cal L}_n^{(j+1)} & = & {\cal L}_n^{(j)} + \gamma_{n(p+1)+j+2}{\cal L}_{n+1}^{(j)} \label{duales1}\\
(z-C) {\cal L}_n^{(j)} & = &  {\cal L}_{n-1}^{(j+1)} +\sum_{s=0}^{p-1} g_{n+s+1,n+1}^{(j)}{\cal L}_{n+s}^{(j+1)}\label{duales2}\\
(z-C){\cal L}_n & = & {\cal L}_{n-1}^{(p)}-\frac{P_{n+1}(C) }{P_n(C)}{\cal L}_{n}^{(p)}\label{duales3}
\end{eqnarray}
\end{lem}

\noindent
\underline{Proof}.- 
In \cite{Dani}, the relation
\begin{equation}
P_{m+1}^{(j)}=P_{m+1}^{(j+1)}+\gamma_{m(p+1)+j+2}P_{m}^{(j+1)},\quad m=-1,0,1,\ldots ,
\label{99}
\end{equation}
was proved (here, $\gamma_{-(p+1)+j+2}=0$). Then
\begin{eqnarray}\label{dual3}
{\cal L}_n^{(j+1)} \left[P_{m+1}^{(j)}\right]&= & {\cal L}_n^{(j+1)} \left[P_{m+1}^{(j+1)}\right]+\gamma_{m(p+1)+j+2} {\cal L}_n^{(j+1)}\left[P_{m}^{(j+1)}\right]\nonumber\\
& = & \left\{
\begin{array}{lll}
0 & , & n\neq m\,,\,m+1\\
1 & , & n=m+1\\
\gamma_{m(p+1)+j+2}& , & n=m\,.
\end{array}
\right.
\end{eqnarray}
Moreover, 
$$
\left({\cal L}_n^{(j)} + \gamma_{n(p+1)+j+2}{\cal L}_{n+1}^{(j)}\right)\left[P_{m+1}^{(j)}\right]=
{\cal L}_n^{(j)}\left[P_{m+1}^{(j)}\right]+ \gamma_{n(p+1)+j+2}{\cal L}_{n+1}^{(j)}\left[P_{m+1}^{(j)}\right]
$$
also drives to (\ref{dual3}). That is, both sides of (\ref{duales1})  coincide on the basis $\{P_{m}^{(j)}\}\,,\, m \in \mathbb{N}\,, $ of the space $\mathcal{P}$ of polynomials. Therefore (\ref{duales1}) is verified.

For $m=0,1,\ldots, $ taking into account (\ref{kernel1}), 
\begin{eqnarray}
\label{duales11}
(z-C) {\cal L}_n^{(j)}\left[  P_{m}^{(j+1)} \right] & = & {\cal L}_n^{(j)}\left[ (z-C) P_{m}^{(j+1)} \right] \nonumber \\
& = & {\cal L}_n^{(j)}\left[ P^{(j)}_{m+1}\right]+\displaystyle\sum_{s=0}^{p-1}g^{(j)}_{m+1,m-s+1}{\cal L}_n^{(j)}\left[P^{(j)}_{m-s}\right]\nonumber \\
& = & \left\{
\begin{array}{lll}
0 &,& n\neq m+1,m,m-1,\ldots, m-p+1\,.\\
1 &,& n=m+1\\
g^{(j)}_{n+s+1,n+1} &,& n=m-s\,,\quad s=0,1,\ldots, p-1\,.
\end{array}
\right.
\end{eqnarray}
On the other hand, 
$$
\left({\cal L}_{n-1}^{(j+1)} +\sum_{s=0}^{p-1} g_{n+s+1,n+1}^{(j)}{\cal L}_{n+s}^{(j+1)}\right)\left[  P_{m}^{(j+1)} \right]=
{\cal L}_{n-1}^{(j+1)} \left[  P_{m}^{(j+1)} \right]+\sum_{s=0}^{p-1} g_{n+s+1,n+1}^{(j)}{\cal L}_{n+s}^{(j+1)}\left[  P_{m}^{(j+1)} \right]
$$
coincides with (\ref{duales11}) for each $m\in \mathbb{N}$. Therefore, (\ref{duales2}) holds. We underline that this is true even if $n=0$ in  (\ref{duales2}), understanding ${\cal L}_{-1}^{(j+1)}=0 $ in this case. 

As in (\ref{duales1})-(\ref{duales2}), we apply both sides of (\ref{duales3}) to a basis of polynomials. Then using (\ref{kernel2}), 
\begin{eqnarray}\label{izquierda}
(z-C){\cal L}_n \left[ P^{(p)}_{m}\right] & = &{\cal L}_n \left[  (z-C)P^{(p)}_{m}\right]\nonumber\\
& = & {\cal L}_n \left[   P_{m+1}- \frac{P_{m+1}(C) }{P_m(C)}P_m\right]\\
& = &  {\cal L}_n \left[   P_{m+1}\right]-\frac{P_{m+1}(C) }{P_m(C)} {\cal L}_n \left[ P_m\right]=
\left\{
\begin{array}{cll}
0 &,& n\neq m,\,m+1\\
1& , & n=m+1\\
-\displaystyle\frac{P_{m+1}(C) }{P_m(C)} & , & n=m.
\end{array}
\right.\nonumber
\end{eqnarray}
Further,
$$
\left({\cal L}_{n-1}^{(p)}-\frac{P_{n+1}(C) }{P_n(C)}{\cal L}_{n}^{(p)}\right)\left[ P^{(p)}_{m}\right]= 
{\cal L}_{n-1}^{(p)}\left[ P^{(p)}_{m}\right]-\frac{P_{n+1}(C) }{P_n(C)}{\cal L}_{n}^{(p)}\left[ P^{(p)}_{m}\right],
$$
which produces exactly the same result as in (\ref{izquierda}). Then (\ref{duales3}) is proved.     $\hfill \square$

In the classic case $p=1$, the functionals of orthogonality $\nu$ and $\nu^{(1)}$, associated respectively with $\{P_n\}$ and the Kernel polynomials $\{P_n^{(1)}\}$, are related by 
$$
\nu^{(1)}=(z-C)\nu\,.
$$
The next result extends this fact to the general case $p\in \mathbb{N}$. In fact, this lemma is equivalently to Theorem~\ref{principal} in the case $j=p$.  

\begin{lem}
\label{lema4}
Let $\nu=\left(\nu_1,\,\nu_2,\,\ldots , \nu_p\right)$ be a vector of $p$-orthogonality for $\{P_n\}$. Then
$$\nu^{(p)}=\left((z-C)\nu_1,\,(z-C)\nu_2,\,\ldots , (z-C)\nu_p\right)$$
is a  vector of $p$-orthogonality for the $p$-Kernel polynomials $\{P^{(p)}_n\}$.
\end{lem}

\noindent
\underline{Proof}.- Due to (\ref{kernel2}), for each $r=1,\,2,\,\ldots,\,p$ we have
\begin{eqnarray*}
\left((z-C)\nu_r\right)\left[ z^kP^{(p)}_{mp+i}\right]&=&\nu_r\left[ z^k(z-C)P^{(p)}_{mp+i}\right]
=\nu_r\left[z^k\left( P_{mp+i+1}- \frac{P_{mp+i+1}(C) }{P_{mp+i}(C)}P_{mp+i}\right)\right]\\
& = &\nu_r\left[ z^k P_{mp+i+1}\right]- \frac{P_{mp+i+1}(C) }{P_{mp+i}(C)}\nu_r\left[z^kP_{mp+i}\right],
\end{eqnarray*}
where
\begin{equation}
\left\{
\begin{array}{rll}
\nu_r\left[ z^k P_{mp+i+1}\right]=0 &,& k\ge 0\,,\quad kp+r\le mp+i+1,\\\\
\nu_r\left[ z^k P_{mp+i}\right]=0 &,& k\ge 0\,,\quad kp+r\le mp+i.
\end{array}
\right.
\label{dos}
\end{equation}
Therefore, if $kp+r\leq mp+i$ then (\ref{dos}) holds and 
$$
\left((z-C)\nu_r\right)\left[z^kP^{(p)}_{mp+i}\right]=0\,,\qquad k\ge 0,\,\quad kp+r\le mp+i\,.
$$
Moreover, using (\ref{kernel1}),
$$
\left((z-C)\nu_r\right)\left[z^kP^{(p)}_{kp+r-1}\right]=\nu_r\left[z^kP_{kp+r}\right]-\frac{P_{kp+r}(C)}{P_{kp+r-1}(C)}\nu_r\left[z^kP_{kp+r-1}\right],
$$
where 
$\nu_r\left[z^kP_{kp+r}\right]=0$ (see (\ref{dos})) and $\nu_r\left[z^kP_{kp+r-1}\right]\neq 0$ (see (\ref{ortogonalidad})). This is,
$$
\left((z-C)\nu_r\right)\left[z^kP^{(p)}_{kp+r-1}\right]=-\frac{P_{kp+r}(C)}{P_{kp+r-1}(C)}\nu_r\left[z^kP_{kp+r-1}\right]\neq 0\,.
$$
This proves that $(z-C)\nu_r$ is the $r$-th entry of a vector of $p$-orthogonality associated with the sequence of polynomials $\{P^{(p)}_n\}$ and, consequently, $\nu^{(p)}$ is one of such vectors. $\hfill \square$

\begin{rem}
\label{rem2}
We underline that, in the case $j=p$, we have proved that the statement of Theorem \ref{principal} is verified independently on the condition \eqref{condition}.
\end{rem}

\begin{lem}\label{teorema2}
For each $j=0,1,\ldots, p$, let $\nu^{(j)}=\left(\nu^{(j)}_1,\,\nu^{(j)}_2\,\ldots , \nu^{(j)}_p\right)$ be a vector of $p$-orthogonality for $\{P^{(j)}_{n}\}$. Then,  for $j=0,1,\ldots, p-1$ we have: 
\begin{itemize}
\item[(a) ]   $\tilde\nu^{(j+1)}=\left(\nu^{(j+1)}_1,\,\ldots , \,\nu^{(j+1)}_{p-1}, (z-C)\nu^{(j)}_1\right)$ is a vector of $p$-orthogonality for $\{P^{(j+1)}_{n}\}$.
\item[(b)]  $\tilde\nu^{(j)}=\left(\nu^{(j)}_1,\nu^{(j+1)}_1\,\ldots , \,\nu^{(j+1)}_{p-1}\right)$ is a vector of $p$-orthogonality for $\{P^{(j)}_{n}\}$.
\end{itemize}
\end{lem}
\noindent
\underline{Proof}.- In the first place, because the first entries of the vector $\tilde\nu^{(j+1)}$ coincide with the corresponding to $\nu^{(j+1)}$, to prove $(a)$ it is enough to check 
\begin{eqnarray}
%\left\{
(z-C)\nu_{1}^{(j)}\left[z^kP_{n}^{(j+1)}\right]& = & 0\,,\quad kp+p\leq n\,,\label{124}\\
(z-C)\nu_{1}^{(j)}\left[z^kP_{(k+1)p-1}^{(j+1)}\right]&\neq & 0\,.\label{125}
%\right.
\end{eqnarray}
Indeed, using (\ref{kernel1}) of Theorem \ref{teorema1},
\begin{equation*}
(z-C)\nu_{1}^{(j)}\left[z^kP_{n}^{(j+1)}\right]=\nu_{1}^{(j)}\left[z^kP_{n+1}^{(j)}\right]+
\sum_{s=0}^{p-1}g^{(j)}_{n+1,n-s+1}\nu_{1}^{(j)}\left[z^kP_{n-s}^{(j)}\right],
\label{123}
\end{equation*}
where 
$
\nu_{1}^{(j)}\left[z^kP_{n+1}^{(j)}\right]=0
$
for $kp+1\leq n+1$ and
$
\nu_{1}^{(j)}\left[z^kP_{n-s}^{(j)}\right]=0
$
for $kp+1\leq n-s\,,\,s=0,1,\ldots , p-1\,.$

Then, 
$$
(z-C)\nu_{1}^{(j)}\left[x^kP_{n}^{(j+1)}\right]=0
$$
for $kp+1\leq n-p+1$ or, what is the same, (\ref{124}) holds.

For a similar reason,
\begin{eqnarray*}
(z-C)\nu_{1}^{(j)}\left[z^kP_{kp+p-1}^{(j+1)}\right]& = &\nu_{1}^{(j)}\left[z^kP_{kp+p}^{(j)}\right]+
\sum_{s=0}^{p-1}g^{(j)}_{kp+p,kp-s+p}\nu_{1}^{(j)}\left[z^kP_{kp-s+p-1}^{(j)}\right]\nonumber \\
& = & g^{(j)}_{kp+p,kp+1}\nu_{1}^{(j)}\left[z^kP_{kp}^{(j)}\right]\neq 0\,,
\end{eqnarray*}
which, taking into account (\ref{asterisco}), gives (\ref{125}). Thus $(a)$ is verified.
 
 In the second place, take $r\in \{1,\ldots ,p-1\}, \,k\in \{0,1,\ldots, \}$ and $n\in \mathbb{N}$. Using (\ref{99}),
\begin{equation}
\nu_r^{(j+1)}\left[z^kP_{n}^{(j)}\right]=\nu_r^{(j+1)}\left[z^kP_{n}^{(j+1)}\right]+
\gamma_{(n-1)(p+1)+j+2}\nu_r^{(j+1)}\left[z^kP_{n-1}^{(j+1)}\right].
\label{100}
\end{equation}
In (\ref{ortogonalidad}) we see
\begin{equation}
\nu_r^{(j+1)}\left[z^kP_{n}^{(j+1)}\right]=\nu_r^{(j+1)}\left[z^kP_{n-1}^{(j+1)}\right]=0\,,\quad  kp+r\leq n-1\,.
\label{otra40}
\end{equation}
Hence (\ref{100}) implies 
\begin{equation}
\nu_r^{(j+1)}\left[z^kP_{n}^{(j)}\right]=0\,,\quad kp+r+1\leq n\,.
\label{101}
\end{equation}
From (\ref{100})-(\ref{otra40}), taking $n=kp+r$,
\begin{equation}
\nu_r^{(j+1)}\left[z^kP_{kp+r}^{(j)}\right]=\gamma_{(kp+r-1)(p+1)+j+2}   \nu_r^{(j+1)}\left[z^kP_{kp+r-1}^{(j+1)}\right]\neq 0\,.
\label{102}
\end{equation}
Since \eqref{ortogonalidad}, we have that (\ref{101}) and (\ref{102}) give 
$
\nu_r^{(j+1)}=\widetilde\nu_{r+1}^{(j)},
$
which is the $(r+1)$-th entry of   vector of $p$-orthogonality for $\{P_{n}^{(j)}\},\,n\in \mathbb{N}$ 
(we recall that $r+1\leq p$).

$\hfill \square$

\section{Proof of Theorem \ref{principal}}\label{section3}

Through Lemma \ref{lema4} and Remark \ref{rem2}, the result is verified  for $j=p$, independent on the factorization \eqref{5}. Then we want to find a Darboux factorization \eqref{5} such that \eqref{vectorj} is a vector of $p$-orthogonality for the corresponding sequence $\{P^{(j)}_n\}$ of polynomials when $j=1,\ldots , p-1$.

We proceed recursively on $j=1,2,\ldots, p$.

\subsection{First step: $j=1$}

In this case, in \eqref{vectorj} we have the vector of functionals
\begin{equation}
\label{(0)}
\nu^{(1)}=\left(\nu_2,\ldots, \nu_p,(z-C)\nu_1\right).
\end{equation}
Due to Lemma \ref{teorema2}, to prove Theorem \ref{principal} it is sufficient to show that $\nu_2,\ldots, \nu_p$ are the first $p-1$ entries of a vector of $p$-orthogonality for $\{P^{(j)}_n\}$, where this sequence of polynomials is corresponding to some Darboux transformation $J^{(1)}$ of $J-CI$. In this step, we  choose $L^{(1)}$ appropriately for our goal. This is, we will see how to fix the entries 
$$
\gamma_2,\,\gamma_{p+3},\ldots, \gamma_{(p-1)p}
$$
of $L^{(1)}$ in \eqref{puntos} with the purpose to define $\mathcal{L}_{m}^{(1)}$ as in \eqref{duales1} and to find 
 $\lambda_{m,k}^{(1)}\in \mathbb{C},\,k=0,1,\ldots,m-1,$ and $\lambda_{m,m-1}^{(1)}\neq 0$ such that
\begin{equation}
\label{(1)}
\nu_{m+1}=\sum_{k=0}^{m-1} \lambda_{m,k}^{(1)}\mathcal{L}_{k}^{(1)}\,,\quad m=1,2,\ldots, p-1
\end{equation}
(see \eqref{escalera}).

Because $\nu$ is a vector of $p$-orthogonality for  $\{P_n\}$, we know
\begin{equation}
\label{(2)}
\nu_{m+1}=\sum_{k=0}^{m} \lambda_{m+1,k}\mathcal{L}_{k}\,,\quad m=1,2,\ldots, p-1.
\end{equation}
For any Darboux factorization, \eqref{duales1} holds taking $j=0$. From this \eqref{(1)} is equivalent to
$$
\nu_{m+1}= \lambda_{m,0}^{(1)}\mathcal{L}_{0}+\sum_{k=1}^{m-1} \left(\gamma_{(k-1)(p+1)+2}\lambda_{m,k-1}^{(1)}+\lambda_{m,k}^{(1)} \right)\mathcal{L}_{k}+\gamma_{(m-1)(p+1)+2}\lambda_{m,m-1}^{(1)}\mathcal{L}_{m}.
$$
Comparing the last expression with \eqref{(2)} for $m=1,\ldots, p-1,$ we have
\begin{eqnarray}
\lambda_{m,0}^{(1)} & = & \lambda_{m+1,0},\label{53} \\
\gamma_{(k-1)(p+1)+2}\lambda_{m,k-1}^{(1)}+\lambda_{m,k}^{(1)} & = & \lambda_{m+1,k}\,,\quad k=1,\ldots, m-1,\label{54}\\
\gamma_{(m-1)(p+1)+2}\lambda_{m,m-1}^{(1)} & = &  \lambda_{m+1,m}.\label{55}
\end{eqnarray}
We rewrite \eqref{53}-\eqref{55} as 
\begin{equation}
\label{57}
L_{m+1}^{(1)}
\left(
\begin{array}{c}
\lambda_{m,0}^{(1)}\\
\lambda_{m,1}^{(1)}\\
\vdots \\
\lambda_{m,m-1}^{(1)}\\
0
\end{array}
\right)=
\left(
\begin{array}{c}
\lambda_{m+1,0}\\
\lambda_{m+1,1}\\
\vdots \\
\lambda_{m+1,m-1}\\
\lambda_{m+1,m}
\end{array}
\right),\,m=1,\ldots, p-1\,.
\end{equation}
In other words, \eqref{53}-\eqref{54} is
\begin{equation*}
\label{otro57}
\left(
\begin{array}{c}
\lambda_{m,0}^{(1)}\\
\lambda_{m,1}^{(1)}\\
\vdots \\
\lambda_{m,m-1}^{(1)}
\end{array}
\right)=
\left(L_{m}^{(1)}\right)^{-1}
\left(
\begin{array}{c}
\lambda_{m+1,0}\\
\lambda_{m+1,1}\\
\vdots \\
\lambda_{m+1,m-1}
\end{array}
\right),\,m=1,\ldots, p-1\,,
\end{equation*}
or, what is the same,
$$
\lambda_{m,k-1}^{(1)} =\lambda_{m+1,k-1}-\gamma_{(k-2)(p+1)+2}\lambda_{m+1,k-2}
+ \gamma_{(k-3)(p+1)+2}\gamma_{(k-2)(p+1)+2}\lambda_{m+1,k-3} -
$$
\begin{equation}
\label{(4)}
\cdots + (-1)^{k-1}\gamma_{2}\gamma_{(p+1)+2}\ldots\gamma_{(k-2)(p+1)+2}\lambda_{m+1,0},\quad 1\leq k\leq m,\quad  1\leq m\leq p-1.
\end{equation}
Furthermore, taking into account the expression of the last row of $\left(L_{m+1}^{(1)}\right)^{-1}$, 
the following relation joint with \eqref{(4)} are equivalent to \eqref{57},
\begin{eqnarray}
\label{(5)}
\lambda_{m+1,m}-\gamma_{(m-1)(p+1)+2}\lambda_{m+1,m-1}
+ \cdots+(-1)^{m}\gamma_{2}\gamma_{(p+1)+2}\ldots\gamma_{(m-1)(p+1)+2}\lambda_{m+1,0}&=&0,\nonumber \\
\quad m=1,\ldots, p-1. &&
\end{eqnarray}
Therefore, the proof of the case $j=1$ is reduced to find entries $\gamma_2,\ldots, \gamma_{(p-1)p}$ of $L^{(1)}$ verifying the condition \eqref{(5)} such that the coefficients $\lambda_{m,k-1}^{(1)}, k=1,\ldots,m,$ provided in \eqref{(4)} define a vector of functionals
$$
\nu^{(1)}=\left(\nu^{(1)}_1,\ldots, \nu^{(1)}_{p-1},(z-C)\nu_1\right)
$$ 
as in \eqref{(0)}.  We proceed recursively for $m=1,2,\ldots, p-1$. 

For $m=1$ in \eqref{(5)} we have 
$
\lambda_{2,1}-\gamma_2\lambda_{2,0}=0\,,
$
where we know from \eqref{condition} that $\Delta_1=\lambda_{2,0}\neq 0$. Then, taking  $\Delta_0:=1$, we can define
\begin{equation*}
%\label{gamma2}
\gamma_2=\lambda_{2,1}\frac{\Delta_0}{\Delta_1}.
\end{equation*}
Moreover, since \eqref{(4)} we define $\lambda^{(1)}_{1,0}=\lambda_{2,0}$ and $\nu^{(1)}_1=\lambda_{2,0}\mathcal{L}_0$ has been constructed. 

Now we will to prove that the first entries of $L^{(1)}$ can be choosen as
\begin{equation}
\gamma_{(m-1)(p+1)+2}=\lambda_{m+1,m}\frac{\Delta_{m-1}}{\Delta_m}\,,\quad m=1,2,\ldots, p-1.
\label{(8)}
\end{equation}
Indeed, \eqref{(8)} is verified for $m=1$. Assume that \eqref{(8)} holds for $m\leq s<p-1$. Assume also that
$$
\gamma_2,\, \gamma_{(p+1)+2},\,\ldots, \gamma_{(s-1)(p+1)+2}
$$
have been choosen verifying \eqref{(5)}. Then $\gamma_{s(p+1)+2}$  can be defined taking $m=s+1$ in \eqref{(5)}, this is,
\begin{eqnarray*}
0 &= &\lambda_{s+2,s+1}-\gamma_{s(p+1)+2}\left[ \lambda_{s+2,s}-\gamma_{(s-1)(p+1)+2}\lambda_{s+2,s-1}+\right.\nonumber \\
& + & \left.\cdots+(-1)^{s+1}\gamma_{2}\gamma_{(p+1)+2}\ldots\gamma_{(s-1)(p+1)+2}\lambda_{s+2,0}\right],
\end{eqnarray*}
where, from \eqref{(8)}, we see that 
$$
\Delta_s\left[\lambda_{s+2,s}-\gamma_{(s-1)(p+1)+2}\lambda_{s+2,s-1}
+\cdots+(-1)^{s+1}\gamma_{2}\gamma_{(p+1)+2}\ldots\gamma_{(s-1)(p+1)+2}\lambda_{s+2,0}\right]
$$
is the development of the determinant $\Delta_{s+1}$  by its last column. Thus 
$$
\gamma_{s(p+1)+2}=\lambda_{s+2,s+1}\frac{\Delta_s}{\Delta_{s+1}}
$$
and \eqref{(8)} is verified in $m=s+1$ and, consequently, for all $m=1,2,\ldots, p-1$. 

In this way the entries 
$$
\gamma_{2},\,\gamma_{(p+1)+2},\,\ldots,\,\gamma_{(p-1)p}
$$
of $L^{(1)}$ are chosen verifying \eqref{(8)} and the coefficients $\lambda_{m,k-1}^{(1)}, k=1,\ldots,m,$ given in \eqref{(4)} define the vector of orthogonality $$
\nu^{(1)}=\left(\nu^{(1)}_1,\ldots, \nu^{(1)}_{p-1},(z-C)\nu_1\right)
$$ 
for a new sequence of polynomials $\{P_n^{(1)}\}$.\\
\subsection{Steps $2,3,\ldots, p-1$.}
 In each one of the following steps, we want to find the first appropriate entries of the corresponding bidiagonal matrix. This is, in the step $j+1$, for $j=1,2,\ldots, p~-~1$, we assume that for each $s=1,2,\ldots ,j$ the entries 
$$
\gamma_{s+1},\,\gamma_{(p+1)+s+1},\,\ldots, \gamma_{(p-s-1)(p+1)+s+1}
$$ 
of $L^{(s)}$ have been defined such that 
$$
\nu^{(s)}=\left( \nu_{s+1},\ldots,\nu_p , (z-C)\nu_1 , \ldots, (z-C)\nu_s\right)
$$
is a vector of $p$-orthogonality for $\{P_n^{(s)}\}$. Then, we want to find the first $p-j-1$ entries
$$
\gamma_{j+2},\,\gamma_{(p+1)+j+2},\,\ldots, \gamma_{(p-j-1)p}
$$ 
of $L^{(j+1)}$ such that $\nu^{(j+1)}$ in \eqref{vectorj} is a vector of $p$-orthogonality for the corresponding sequence $\{P_n^{(j+1)}\}$.
Due to the case $j+1=p$ is solved in Lemma \ref{lema4}, we assume $j\in \{1,2,\ldots, p-2\}$ in the following. 
We differentiate two kind of entries in $\nu^{(j+1)}$. This is, we denote 
$$
\nu^{(j+1)}=\left( \nu^{(j+1)}_1,\ldots,\nu^{(j+1)}_p  \right)
$$
where 
\begin{equation}
\label{otro61}
\nu^{(j+1)}_k=\left\{
\begin{array}{lcl}
\nu_{j+k+1} &,& k=1,\ldots, p-j-1,\\\\
(z-C)\nu_{j+k+1-p} &,& k=p-j, \ldots, p.
\end{array}
\right.
\end{equation}

In the first place, we analyze the entries $\nu^{(j+1)}_k\,,\,k=1,\ldots,p-j-1$, for which we want to define 
$\mathcal L_{s}^{(j+1)},\,s=0,\ldots, k-1,$ as in \eqref{duales1} and to find 
$
\lambda_{k,0}^{(j+1)}, \ldots, \lambda_{k,k-1}^{(j+1)}
$
such that 
\begin{equation}
\label{111}
\nu_k^{(j+1)}=\sum_{s=0}^{k-1}\lambda_{k,s}^{(j+1)}\mathcal L_{s}^{(j+1)}.
\end{equation}
Since \eqref{vectorj} and \eqref{otro61} we know that
$$
\nu_k^{(j+1)}=\nu_{k+1}^{(j)}\,,\quad k=1,\ldots, p-j-1 \,.
$$
Hence,
\begin{equation}
\label{222}
\nu_k^{(j+1)}=\sum_{r=0}^{k}\lambda_{k+1,r}^{(j)}\mathcal L_{r}^{(j)}\,,\quad k=1,\ldots, p-j-1 \,.
\end{equation}
Using \eqref{duales1} in \eqref{111},
\begin{eqnarray*}
\nu_k^{(j+1)} &= & \sum_{s=0}^{k-1}\lambda_{k,s}^{(j+1)}\left(\mathcal L_{s}^{(j)}+\gamma_{s(p+1)+j+2}\mathcal L_{s+1}^{(j)} \right)=\\
 & = &  \sum_{s=0}^{k-1}\lambda_{k,s}^{(j+1)}\mathcal L_{s}^{(j)}+\sum_{s=1}^{k}\gamma_{(s-1)(p+1)+j+2}\lambda_{k,s-1}^{(j+1)}\mathcal L_{s}^{(j)}=\\
  =\lambda_{k,0}^{(j+1)}\mathcal L_{0}^{(j)}&+& \sum_{s=1}^{k-1}\left(  \lambda_{k,s}^{(j+1)}+\gamma_{(s-1)(p+1)+j+2}\lambda_{k,s-1}^{(j+1)} \right)\mathcal L_{s}^{(j)}+ \gamma_{(k-1)(p+1)+j+2}\lambda_{k,k-1}^{(j+1)} \mathcal L_{k}^{(j)}.
\end{eqnarray*}
Comparing with \eqref{222},
\begin{equation}
\label{333}
\lambda_{k+1,s}^{(j)}=
\left\{
\begin{array}{rlll}
& \lambda_{k,0}^{(j+1)}\,, & s=0, \\\\
 &  \lambda_{k,s}^{(j+1)}+  \gamma_{(s-1)(p+1)+j+2}\lambda_{k,s-1}^{(j+1)}\,,&  s=1,\ldots ,k-1,\\\\
 & \gamma_{(k-1)(p+1)+j+2}\lambda_{k,k-1}^{(j+1)}\,, & s=k\,.
\end{array}
\right.
\end{equation}
Similarly to what was done in \eqref{53}-\eqref{55},
\eqref{333} can be rewritten as

\begin{equation}
\label{otro3}
L_{k+1}^{(j+1)}
\left(
\begin{array}{c}
\lambda_{k,0}^{(j+1)}\\
\lambda_{k,1}^{(j+1)}\\
\vdots \\
\lambda_{k,k-1}^{(j+1)}\\
0
\end{array}
\right)=
\left(
\begin{array}{c}
\lambda_{k+1,0}^{(j)}\\
\lambda_{k+1,1}^{(j)}\\
\vdots \\
\lambda_{k+1,k-1}^{(j)}\\
\lambda_{k+1,k}^{(j)}
\end{array}
\right),\,k=1,\ldots, p-j-1\,,
\end{equation}
or, what is the same,
\begin{equation}
\label{otro4}
L_{k}^{(j+1)}
\left(
\begin{array}{c}
\lambda_{k,0}^{(j+1)}\\
\lambda_{k,1}^{(j+1)}\\
\vdots \\
\lambda_{k,k-1}^{(j+1)}
\end{array}
\right)=
\left(
\begin{array}{c}
\lambda_{k+1,0}^{(j)}\\
\lambda_{k+1,1}^{(j)}\\
\vdots \\
\lambda_{k+1,k-1}^{(j)}
\end{array}
\right),\,k=1,\ldots, p-j-1\,,
\end{equation}
with the additional condition
\begin{equation}
\label{otro5}
\lambda_{k+1,k}^{(j)} = \gamma_{(k-1)(p+1)+j+2}\lambda_{k,k-1}^{(j+1)}\,,\quad k=1,\ldots,p-j-1\,.
\end{equation}
Therefore, finding
$$
\lambda_{k,0}^{(j+1)},\ldots, \lambda_{k,k-1}^{(j+1)}
$$
as in \eqref{111} comes down to choose the entries 
$$
 \gamma_{j+2},\ldots,  \gamma_{(p-j-1)p}
$$
of $L_{p-j}^{(j+1)}$ with the aim of \eqref{otro4}-\eqref{otro5} take place. We note that \eqref{otro4} defines the coefficients 
$\lambda_{k,s}^{(j+1)},\,s=0,\ldots, k-1,$ because $L_{p-j}^{(j+1)}$  is an invertible matrix. Then, \eqref{otro5} is equivalent to the fact that the last row of 
$
\left( L_{k+1}^{(j+1)} \right)^{-1}$ multiplied by   
$
\left(
\lambda_{k+1,0}^{(j)},\ldots, \lambda_{k+1,k}^{(j)}
\right)
$
vanishes. This is, 
\begin{eqnarray}
\label{otro6}
\lambda_{k+1,k}^{(j)}&- &
\lambda_{k+1,k-1}^{(j)}\gamma_{(k-1)(p+1)+j+2}+
\lambda_{k+1,k-2}^{(j)}\gamma_{(k-2)(p+1)+j+2}\gamma_{(k-1)(p+1)+j+2}\\
-\cdots&+&
(-1)^k\lambda_{k+1,0}^{(j)}\gamma_{j+2}\gamma_{(p+1)+j+2}\ldots\gamma_{(k-1)(p+1)+j+2}=0\,,\quad k=1,2,\ldots, p-j\,.\nonumber
\end{eqnarray}

Now we show that in the above conditions the following matrix equality is verified for $j=1,2,\ldots, p-1$ and $s=1,2,\ldots$,
\begin{eqnarray} 
\label{otro0}  
L^{(1)}_s\ldots L^{(j)}_s
\left(
\begin{array}{cccccccc}
\lambda^{(j)}_{2,0} & \lambda^{(j)}_{3,0}  &\cdots  &\cdots  & \cdots & \lambda^{(j)}_{s+1,0}\\
\lambda^{(j)}_{2,1} &\lambda^{(j)}_{3,1}  & \cdots  & \cdots  &\cdots &\lambda^{(j)}_{s+1,1}\\
0 & \lambda^{(j)}_{3,2} &\ddots& && \vdots\\
\vdots & 0 &\ddots &\ddots &\\
\vdots & \vdots & \ddots &&&\vdots  \\
0 & 0 & \cdots &0 & \lambda^{(j)}_{s,s-1} & \lambda^{(j)}_{s+1, s-1}\\
\end{array}
\right)\nonumber \\\nonumber \\
=\left(
\begin{array}{cccccccc}
\lambda_{j+2,0} & \lambda_{j+3,0}  &\cdots  &\cdots  & \lambda_{s,0}  &\cdots & \lambda_{j+s+1,0}\\
\vdots & \vdots & & & \vdots & \\
\lambda_{j+2,j+1} &\lambda_{j+3,j+1}  & \cdots  & \cdots  &\lambda_{s,j+1} &\cdots &\lambda_{j+s+1,j+1}\\
0 & \lambda_{j+3,j+2} &\ddots&& \vdots&& \vdots\\
\vdots & 0 &\ddots &\ddots && \cdots \\
\vdots & \vdots & \ddots &\ddots & \vdots&& \vdots  \\
0 & 0 & \cdots &0 & \lambda_{s,s-1} & \cdots & \lambda_{j+s+1, s-1}\\
\end{array}
\right).
\end{eqnarray}
In fact, as in \eqref{otro3}, it is easy to see
$$
L_{k+1}^{(r)}
\left(
\begin{array}{c}
\lambda_{k,0}^{(r)}\\
\lambda_{k,1}^{(r)}\\
\vdots \\
\lambda_{k,k-1}^{(r)}\\
0
\end{array}
\right)=
\left(
\begin{array}{c}
\lambda_{k+1,0}^{(r-1)}\\
\lambda_{k+1,1}^{(r-1)}\\
\vdots \\
\lambda_{k+1,k-1}^{(r-1)}\\
\lambda_{k+1,k}^{(r-1)}
\end{array}
\right),\,k=1,\ldots, p-r,\quad r=1,2,\ldots, j+1\,.
$$
Then, considering each infinite matrix $L^{(r)}$,
\begin{equation}
\label{1111}
L^{(r)}
\left(
\begin{array}{c}
\lambda_{k,0}^{(r)}\\
\lambda_{k,1}^{(r)}\\
\vdots \\
\lambda_{k,k-1}^{(r)}\\
0\\
0\\
\vdots 
\end{array}
\right)=
\left(
\begin{array}{c}
\lambda_{k+1,0}^{(r-1)}\\
\lambda_{k+1,1}^{(r-1)}\\
\vdots \\
\lambda_{k+1,k-1}^{(r-1)}\\
\lambda_{k+1,k}^{(r-1)}\\
0\\
\vdots 
\end{array}
\right),\,k=1,\ldots, p-r,\quad r=1,2,\ldots, j+1\,.
\end{equation}
Applying iteratively $L^{(r)},\ldots , L^{(1)}$ to \eqref{1111} and then taking $r=j$ we arrive to 
$$
L^{(1)}\cdots L^{(j)}
\left(
\begin{array}{c}
\lambda_{k,0}^{(j)}\\
\vdots \\
\lambda_{k,k-1}^{(rj)}\\
0\\
\vdots 
\end{array}
\right)=
\left(
\begin{array}{c}
\lambda_{j+k,0}\\
\vdots \\
\vdots \\
\lambda_{j+k,j+k-1}\\
0\\
\vdots 
\end{array}
\right),\,k=1,\ldots, p-j.
$$
From this, for $s\in \mathbb{N}$, taking $k=2,3,\ldots, s+1$, 
$$
L^{(1)}\cdots L^{(j)}
\left(
\begin{array}{ccccc}
 \lambda_{2,0}^{(j)}&\lambda_{3,0}^{(j)}&\cdots & \lambda_{s+1,0}^{(j)}\\
 \lambda_{2,1}^{(j)}&\lambda_{3,1}^{(j)}&\cdots & \lambda_{s+1,1}^{(j)}\\
0&\lambda_{3,2}^{(j)}& \ddots  & \vdots \\
\vdots & 0 &  \cdots & \vdots \\
 & \vdots & \ddots & \lambda_{s+1,s}^{(j)}\\
& & & 0\\
& & & \vdots 
\end{array}
\right)=
\left(
\begin{array}{ccccc}
 \lambda_{j+2,0}&\lambda_{j+3,0}&\cdots & \lambda_{j+s+1,0}\\
 \vdots & \vdots & & \vdots \\
 \lambda_{j+2,j+1}&\lambda_{j+3,j+1}&\cdots & \lambda_{j+s+1,j+1}\\
0&\lambda_{j+3,j+2}& \ddots  & \vdots \\
\vdots & 0 &  \cdots & \vdots \\
 & \vdots & \ddots & \lambda_{j+s+1,j+s}\\
& & & 0\\
& & & \vdots 
\end{array}
\right).
$$
Therefore, due to 
$
L^{(1)}\cdots L^{(j)}
$
is an infinite lower triangular matrix, we have 
$$
\left(
L^{(1)}\cdots L^{(j)}
\right)_s=
L_s^{(1)}\cdots L_s^{(j)}
$$
and we arrive to \eqref{otro0}.

As a consequence of \eqref{otro0},
$$
\Delta^{(j)}_s =\left|
\begin{array}{ccccccc}
\lambda^{(j)}_{2,0} & \lambda^{(j)}_{3,0}  &\cdots  &\cdots  & \cdots & \lambda^{(j)}_{s+1,0}\\
 \lambda^{(j)}_{2,1} &\lambda^{(j)}_{3,1}  & \cdots  & \cdots  &\cdots &\lambda^{(j)}_{s+1,1}\\
0 & \lambda^{(j)}_{3,2} &\ddots& && \vdots\\
\vdots & 0 &\ddots &\ddots &\\
\vdots & \vdots & \ddots &&&\vdots  \\
0 & 0 & \cdots &0 & \lambda^{(j)}_{s,s-1} & \lambda^{(j)}_{s+1, s-1}\\
\end{array}
\right| \,,\quad j=1,\ldots,p-1,\quad s\in \mathbb{N}.
$$

Taking $k=1$ in \eqref{otro6}, we have 
$
\lambda_{2,1}^{(j)}-\lambda_{2,0}^{(j)}\gamma_{j+2}=0
$.
Because we know that $\lambda_{2,1}^{(j)}\neq 0$ (see \eqref{escalera}) and 
$\Delta_1^{(j)}=\lambda_{2,0}^{(j)}\neq 0$, it is possible to take 
$$
\gamma_{j+2}=\lambda_{2,1}^{(j)}\frac{\Delta_0^{(j)}}{\Delta_1^{(j)}}
$$
(where we define $\Delta_0^{(j)}:=1$). Iterating the procedure, assuming 
\begin{equation}
\label{otro7}
\gamma_{m(p+1)+j+2}=\lambda_{m+2,m+1}^{(j)}\frac{\Delta_m^{(j)}}{\Delta_{m+1}^{(j)}}\,,\quad m=0,1,\ldots ,s,
\end{equation}
with $s<p-j-2$ and taking $k=s+2$ in \eqref{otro6}, 
\begin{multline*} 
\lambda_{s+3,s+2}^{(j)}-\gamma_{(s+1)(p+1)+j+2}\Big[
\lambda_{s+3,s+1}^{(j)}-\lambda_{s+3,s}^{(j)} \gamma_{s(p+1)+j+2}+\\ \cdots+
(-1)^{s+1}\lambda_{s+3,0}^{(j)}\gamma_{j+2}\ldots \gamma_{s(p+1)+j+2}\Big]=0.
\end{multline*}
where
$
\gamma_{(s+1)(p+1)+j+2}
$
can be defined as 
\begin{equation}
\label{otro8}
\gamma_{(s+1)(p+1)+j+2}=\lambda_{s+3,s+2}^{(j)}\frac{\Delta_{s+1}^{(j)}}{
\Delta_{s+1}^{(j)}
\left[
\lambda_{s+3,s+1}^{(j)}-\cdots+
(-1)^{s+1}\lambda_{s+3,0}^{(j)}\gamma_{j+2}\ldots \gamma_{s(p+1)+j+2}\right]}.
\end{equation}
Further, from \eqref{otro7} it is easy to check that the denominator in \eqref{otro8} is the development of the determinant $\Delta_{s+2}^{(j)}$ by its last column. Thus $\gamma_{m(p+1)+j+2}$ can be defined as in \eqref{otro7} for $m=0,1,\ldots ,p-j-2$. 

Finally, we study the entries $\nu^{(j+1)}_k\,,\,k=p-j,\ldots,p$ of $\nu^{(j+1)}$. We want to prove 
\begin{eqnarray}
(z-C)\nu_{j+k+1-p}\left[z^sP_n^{(j+1)}  \right]=0,&&n\ge sp+k\label{(otro1)}\\
\nonumber\\
(z-C)\nu_{j+k+1-p}\left[z^sP_{sp+k-1}^{(j+1)}  \right]\neq 0, \label{(otro11)}
\end{eqnarray}
because this means, from \eqref{ortogonalidad}, that $(z-C)\nu_{j+k+1-p}$ is the entry $k$ of a vector of $p$-orthogonality  for $\{P_n^{(j+1)}\}$.

Lemma \ref{lema4} implies 
\begin{equation}
\left.
\begin{array}{rcl}
(z-C)\nu_{j+k+1-p}\left[z^sP_m^{(p)}  \right]=0,&&m\ge sp+j+k+1-p\\\\
(z-C)\nu_{j+k+1-p}\left[z^sP_{sp+j+k-p}^{(p)}  \right]\neq 0.
\label{(otro2)}
\end{array}
\right\}
\end{equation}
In addition, since \eqref{99}, for any Darboux factorization, 
$$
\left(
\begin{array}{c}
P_0^{(j+1)}\\
P_1^{(j+1)}\\
\vdots
\end{array}
\right)=
L^{(j+2)}
\left(
\begin{array}{c}
P_0^{(j+2)}\\
P_1^{(j+2)}\\
\vdots
\end{array}
\right)=\cdots =L^{(j+2)}L^{(j+3)}\ldots L^{(p)}
\left(
\begin{array}{c}
P_0^{(p)}\\
P_1^{(p)}\\
\vdots
\end{array}
\right).
$$
Thus, $P_n^{(j+1)}$ can be expressed in terms of the entries in the $(n+1)$-th row of $L^{(j+2)}L^{(j+3)}\ldots L^{(p)}$ and the sequence $\{P^{(p)}_n\}$. From this, taking into account that $L^{(j+2)}L^{(j+3)}\ldots L^{(p)}$ is a lower triangular $(p-j)$-banded matrix, we see
$$
P^{(j+1)}_n=\sum_{r=n+j-p+1}^n \alpha_rP^{(p)}_r,\quad n\ge 0 \quad (\text{with }\alpha_n=1). 
$$
Hence,
\begin{equation}
\label{(otro3)}
(z-C)\nu_{j+k+1-p}\left[z^sP_n^{(j+1)} \right]=\sum_{r=n+j-p+1}^n \alpha_r(z-C)\nu_{j+k+1-p}\left[ z^sP_r^{(p)} \right]
\end{equation}
and, using \eqref{(otro2)}, we see that each term on the right hand side of \eqref{(otro3)} vanishes when $n\geq sp+k$. Thus \eqref{(otro1)} is verified. 

To see \eqref{(otro11)}, if $n=sp+k-1$ in \eqref{(otro3)}, then using \eqref{(otro2)} on the right hand side of \eqref{(otro3)} we have 
$$
\sum_{r=sp+k+j-p}^{sp+k-1} \alpha_r(z-C)\nu_{j+k+1-p}\left[ z^sP_r^{(p)} \right]=
\alpha_{sp+k+j-p}(z-C)\nu_{j+k+1-p}\left[ z^sP_{sp+j+k-p}^{(p)} \right]\neq 0\,.
$$
 \hfill $\square$

\begin{rem}
Note that, if $\Delta_s^{(j)}\neq 0$ for $s=0,1,\ldots,p-m-1$ and $j=1,\cdots, m-1,$ then there exist $m$ bidiagonal matrices $L^{(1)}, \ldots, L^{(m)}$ such that 
\begin{equation*}
J-CI=L^{(1)}\cdots L^{(m)}\tilde LU,
\end{equation*}
where $\tilde L$ is a lower triangular matrix (non bidiagonal, in general). However, if 
$\Delta_s^{(m)}= 0$ for some $s\in\{0,\ldots,p-m-2\}$  then we can not assure the existence of a Darboux factorization \eqref{5} of $J-CI$ such that $\nu^{(m+1)}$ is defined as in \eqref{vectorj}.
\end{rem}

\end{document}